# DISCUSSION OF "BREAKDOWN AND GROUPS" BY P. L. DAVIES AND U. GATHER

By Marc G. Genton[1] and André Lucas

*Texas A&M University and Vrije Universiteit*

**1. Introduction.** In their interesting paper Davies and Gather draw our attention to what they call the "small print" in definitions of breakdown. Working from a formal group structure and a notion of equivariance, they show by a number of examples that a definition of breakdown may be void if not accompanied by a reasonable and precise group structure. This leads them to what we would label their key remark in Section 6: "We know of no situation not based on equivariance considerations where it can be shown that the highest breakdown point for a class of reasonable functionals is less than 1."

Though we agree with their general point that one has to take care not to come up with void definitions, or put differently, to make the small print explicit, we want to draw attention to the relation of their results to an alternative definition of breakdown. In particular, we claim that a different perspective on the notion of breakdown may resolve some of the small print issues.

The definition used by Davies and Gather in their equation (2.4) is a standard one and has its roots in the domain of location and scale estimation. As we argued in Genton and Lucas (2003), it is less useful in a setting with dependent data. For example, in a simple autoregression (AR) of order 1,

$$(1) \qquad Y_t = \theta Y_{t-1} + e_t, \qquad \theta \in (-1,1), e_t \sim N(0,1),$$

the ordinary least squares (OLS) estimator for $\theta$ is driven to zero by replacing one of the $Y_t$'s by an arbitrarily large number. Note that the OLS estimator thus tends to the center rather than the edge of the parameter space. Still, most people would agree that the estimator has lost its usefulness if only one extreme outlier is added. The reason is that the estimator no longer conveys useful information on the uncontaminated data. It is this latter notion that we want to put to the fore.

---

[1]Supported in part by NSF Grant DMS-02-04297.

**2. Breakdown point for (in)dependent observations.** First, we would like to acknowledge that the breakdown definition as introduced in Genton and Lucas (2003) is subject to criticism raised by Davies and Gather (personal communication). One can easily construct an example of an estimator with breakdown point of 1 that would lose its information content on the uncontaminated process upon the addition of only one outlier. This is mainly due to the lack of a limiting operation in our original definition. Therefore, for the sake of this comment we introduce the following slightly adapted and simplified version of the definition in Genton and Lucas (2003).

Let $Y$ denote a vector containing the sample of observations, and let $\mathcal{Y}$ denote the set of allowable samples. For example, in the asymptotic case $Y$ might be a specific AR(1) process, while $\mathcal{Y}$ is the set of all stationary AR(1) processes. In a finite sample, $Y$ might be a specific vector in $\mathbb{R}^n$, while $\mathcal{Y}$ is equal to $\mathbb{R}^n$. Let $Z_k^\zeta$ be an additive outlier process consisting of $k$ outliers of magnitude $\zeta$, such that we observe $Y + Z_k^\zeta$ rather than $Y$. To formalize the notion of information content on the uncontaminated process, we introduce the concept of badness set, which in this case we define as

$$(2) \qquad R^*(Z_k^\zeta, \mathcal{Y}) = \{\theta(Y + Z_k^\zeta) | Y \in \mathcal{Y}\},$$

where $\theta(\cdot)$ denotes the Fisher consistent estimator functional. Let $\mu$ denote an appropriate measure for the badness set. In most cases the Lebesgue measure suffices. Then we define the breakdown point of an estimator as

$$(3) \qquad \mathrm{bdp} = \frac{1}{n} \min\left\{ k - 1 \Big| \text{for all compact } \mathcal{Y}' \subset \mathcal{Y}: \inf_{Z_k^\zeta} \mu(R^*(Z_k^\zeta, \mathcal{Y}') \cap R^*(0, \mathcal{Y}')) = 0 \right\}.$$

An extension to the asymptotic case is straightforward. To see how the definition works, consider the regression example in Section 6 of Davies and Gather. We have $\mathcal{Y} = \mathbb{R}^{n \times 2}$ and $R^*(0, \mathcal{Y}) = [-n, n]$. The estimator is given by $\theta(Y) = \max(-n, \min(n, \theta^{\mathrm{OLS}}(Y)))$, with $\theta^{\mathrm{OLS}}(Y)$ the standard OLS estimator. We set $\mu$ to the standard Lebesgue measure. By taking $k = 1$ and letting the size of the outlier ($\zeta$) diverge, the intersection of the two badness sets in the definition becomes $\{n\}$ or $\{-n\}$, which is a singleton with Lebesgue measure zero. Therefore, the estimator has broken according to our new definition. This appears reasonable as the estimator no longer conveys information about possibly uncontaminated samples.

**3. Time series.** The advantages of a different perspective on breakdown become even more apparent in the time series setting. Again consider our AR(1) example from (1). In the asymptotic case, define the i.i.d. additive outlier process $Z_{p,t}^\zeta$ with $\mathrm{P}[Z_{p,t}^\zeta = \zeta] = \mathrm{P}[Z_{p,t}^\zeta = -\zeta] = p/2$, and $Z_{p,t}^\zeta = 0$ otherwise. Figure 1 presents plots of the badness set $R^*(Z_p^\zeta, \mathcal{Y})$ associated with three estimators of $\theta$ as a function of $\zeta$ for $p = 5\%, 25\%, 50\%$. Here $\mathcal{Y}$ is the set of all stationary AR(1) processes; see the comment in the discussion below. We set $\mu$ to the standard Lebesgue measure.

The first estimator is the OLS estimator which in the above setting yields the badness set (2) based on the explicit expression

$$\theta_{\mathrm{OLS}}(Y + Z_p^\zeta) = \frac{\theta}{1 + p(1-\theta^2)\zeta^2}. \tag{4}$$

Letting the size of the outliers ($\zeta$) diverge, we see that unless $p = 0$, the estimator $\theta_{\mathrm{OLS}}$ tends to zero and the corresponding badness set becomes $\{0\}$; see the first row of Figure 1. Therefore, the asymptotic breakdown point of the OLS estimator for the AR(1) parameter $\theta$ is 0 in the setting described above.

The second estimator is the least median of squares (LMS) estimator of $\theta$. It yields a badness set (2) based on the expression $\theta_{\mathrm{LMS}}(Y + Z_p^\zeta) = \arg\min_{\tilde\theta \in [-1,1]} c$ under the constraint

$$\frac{1}{2} = (1-p)^2 \chi^2\left(\frac{c}{\tau^2}; 0\right) + p(1-p)\left[\chi^2\left(\frac{c}{\tau^2}; \frac{1}{\tau^2}\zeta^2\right) + \chi^2\left(\frac{c}{\tau^2}; \frac{\tilde\theta^2}{\tau^2}\zeta^2\right)\right]$$

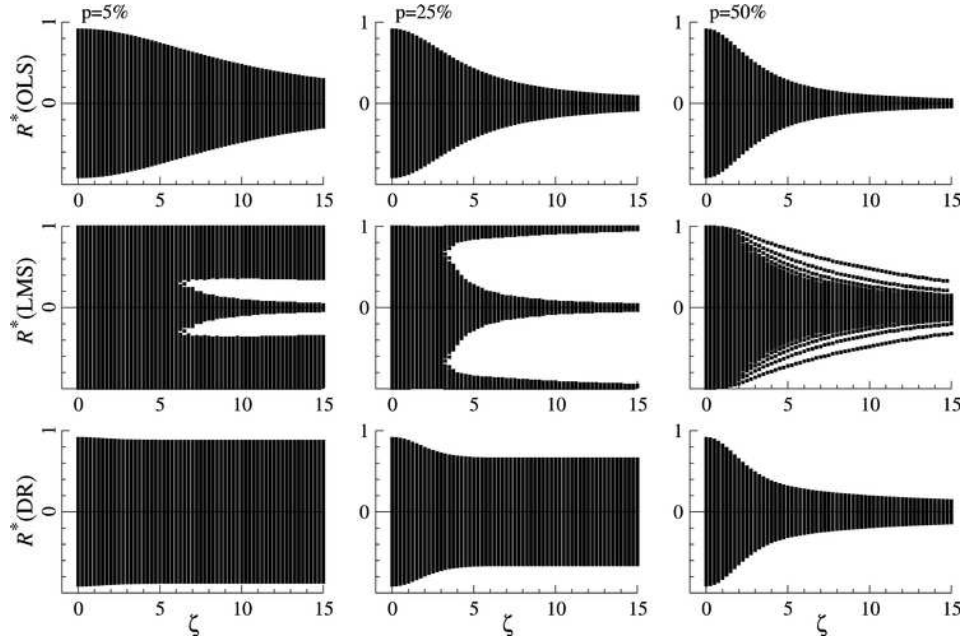

Fig. 1. *Plots of the badness set $R^*(Z_p^\zeta, \mathcal{Y})$ associated with three estimators of $\theta$ in the AR(1) as a function of $\zeta$ for $p = 5\%$, $25\%$, $50\%$: OLS* (top); *LMS* (middle); *DR* (bottom).

(5)
$$+ \frac{p^2}{2}\left[\chi^2\left(\frac{c}{\tau^2};\frac{(1-\tilde{\theta})^2}{\tau^2}\zeta^2\right) + \chi^2\left(\frac{c}{\tau^2};\frac{(1+\tilde{\theta})^2}{\tau^2}\zeta^2\right)\right],$$

where $\tau^2 = 1 + (\theta - \tilde{\theta})^2/(1-\theta^2)$ and $\chi^2(x;\delta^2)$ denotes the cumulative distribution function evaluated at $x$ of a chi-square random variable with noncentrality parameter $\delta^2$. The second row of Figure 1 indicates that the badness set for $p = 5\%$ still takes a continuum of values, whereas it tends to the set $\{-1, 0, +1\}$ for $p = 25\%$. For $p = 50\%$, the badness set collapses to $\{0\}$ as $\zeta$ diverges. Therefore, letting the size of the outliers ($\zeta$) diverge, the asymptotic breakdown point of the LMS estimator for the AR(1) parameter $\theta$ can be computed from (5) to be 22.1% in the setting described above.

The third estimator is the deepest regression (DR) estimator of $\theta$ defined by $\text{median}_t(Y_t/Y_{t-1})$. Under the additive outlier process described above, we need to consider the distribution of $(Y_t + Z^\zeta_{p,t})/(Y_{t-1} + Z^\zeta_{p,t-1})$. It yields a badness set (2) based on the expression $\theta_{\text{DR}}(Y + Z^\zeta_p)$ given by the value $c$ satisfying

$$\frac{1}{2} = (1-p)^2 G(c;0,0)$$

(6)
$$+ \frac{p(1-p)}{2}[G(c;\zeta,0) + G(c;0,\zeta) + G(c;-\zeta,0) + G(c;0,-\zeta)]$$

$$+ \frac{p^2}{4}[G(c;\zeta,\zeta) + G(c;\zeta,-\zeta) + G(c;-\zeta,\zeta) + G(c;-\zeta,-\zeta)],$$

where $G(x;a,b)$ is the cumulative distribution function evaluated at $x$ of the ratio of two correlated normal random variables with means $a$ and $b$, variances $1/(1-\theta^2)$ and correlation $\theta$ [see Hinkley (1969)]. The third row of Figure 1 indicates that the badness set still takes a continuum of values for $p = 5\%$ and $p = 25\%$, whereas it collapses to $\{0\}$ for $p = 50\%$ as $\zeta$ diverges. Thus, the asymptotic breakdown point of the DR estimator for the AR(1) parameter $\theta$ can be computed from (6) to be 50% in the setting described above.

It is interesting to note that the breakdown points of the LMS and DR estimators are markedly different for the AR(1) process above, whereas they are the same (50%) in the setting of simple regression. This indicates that our definition of breakdown allows us to distinguish between various robust estimators in the time series setting.

**4. Discussion.** The definition in (3) appears less dependent on a group structure than the definition used by Davies and Gather. Of course, also the definition in (3) has its limitations. For example, the definition cannot be used if one wants to assess the breakdown of an estimator at a specific sample, that is, if $\mathcal{Y}$ is a singleton. The main drawback of conditioning on the sample is that one has to be very explicit about the region toward which the estimator breaks down, for example, to the edge of the parameter space. This may not be trivial for dependent data, as was shown in the AR(1) example for LMS. Moreover, conditioning the breakdown behavior on a specific sample may relate more to properties of the sample rather than of the estimator. The breakdown notion in (3) based on information revelation about the possible uncontaminated samples resolves this issue. That notion, however, can most easily be operationalized if there is a continuum of possible samples, which suffices for most cases studied in the literature.

A second possible limitation of (3) is that the user has to be explicit about the set $\mathcal{Y}$ of possible samples (or processes) $Y$. For example, if we consider stationary AR(1) processes in the asymptotic setting, the (asymptotic) breakdown point of the OLS estimator is 0. If, however, we consider AR(1) processes characterized by $\theta \in [-1, 1]$, the breakdown point is 1: the OLS estimator retains information about the distinction between stationary processes and processes with $\theta$ arbitrarily close to 1. In that sense the estimator does not break down, while it has broken down if one only wants to distinguish between alternative stationary processes; see Figure 1.

Finally, the definition in (3) is not very explicit about the measure $\mu$. As mentioned, the Lebesgue measure suffices in most cases of practical interest. Despite the fact that empirical data have finite precision, one can work under the assumption that $Y$ lies in a continuum to derive the breakdown properties of the estimator. The properties derived are usually also relevant for a setting with finite precision data. We do not exclude, however, that examples can be constructed where the Lebesgue measure is inappropriate. For example, the parameter space may be discrete and finite. In such cases, alternative measures $\mu$ must be used. Additionally, the restriction that the measure of the intersection of badness sets is zero may have to be replaced by something more complicated, like an infimum of $\inf_{Z_k^\zeta} \mu(R^*(Z_k^\zeta, \mathcal{Y}') \cap R^*(0, \mathcal{Y}'))$ over $k$.

The ideas and cautionary remarks in the paper of Davies and Gather are important and relevant. Effectively, they promote that breakdown is only a useful notion for "sensible" estimators and argue that equivariance is the crucial notion here. We argued that they mainly build on a restricted notion of breakdown. The focus of future research should be put on developing alternative definitions of breakdown that are less susceptible to the criticisms raised by Davies and Gather. The definition in (3) is such an attempt and tries to formalize the phenomena illustrated in Figure 1. In finite samples it is still susceptible to counterexamples, for example, $\theta(Y) = \max(-n, \min(n, \theta^{\mathrm{OLS}}(Y))) + 2(\mathrm{frac}(Y_1) - 1)/n$, where $\mathrm{frac}(x)$ denotes the fractional part of $x$, for Davies and Gather's example in Section **??**, but the examples become increasingly contrived. Moreover, in the asymptotic setting the small print issue appears to become even smaller, especially if we limit ourselves to estimators that are consistent and satisfy some form of continuity in the observations. Further developments along these lines appear promising.

DEPARTMENT OF STATISTICS
TEXAS A&M UNIVERSITY
COLLEGE STATION, TEXAS 77843-3143
USA
E-MAIL: genton@stat.tamu.edu

DEPARTMENT OF FINANCE
ECO/FIN, VRIJE UNIVERSITEIT
DE BOELELAAN 1105
1081HV AMSTERDAM
THE NETHERLANDS
E-MAIL: alucas@feweb.vu.nl